\documentclass[12pt]{amsart}

\usepackage{amsthm, amsmath, amssymb, mathrsfs}
\usepackage{enumerate, microtype}
\usepackage{cite}

\usepackage[colorlinks]{hyperref}

\makeatletter
\@namedef{subjclassname@2010}{%
	\textup{2010} Mathematics Subject Classification}
\makeatother

\allowdisplaybreaks

\theoremstyle{plain}
\newtheorem{thm}{Theorem}[section]

\newtheorem{lem}[thm]{Lemma}
\newtheorem{prop}[thm]{Proposition}
\theoremstyle{definition}

\theoremstyle{remark}
\newtheorem{rem}[thm]{Remark}
\numberwithin{equation}{section}

\title[Topological Center of the Double Dual of $A_\Phi(G)$]{Topological Center of the Double Dual of the Orlicz Fig\`{a}-Talamanca Herz Algebra}
\author[A. Dabra]{Arvish Dabra$^\star$}
\address{Arvish Dabra,\newline\indent Department of Mathematics,\newline\indent Indian Institute of Technology Delhi,\newline\indent New Delhi - 110016, India.}
\email{arvishdabra3@gmail.com}
\author[N. S. Kumar]{N. Shravan Kumar}
\address{N. Shravan Kumar,\newline\indent Department of Mathematics,\newline\indent Indian Institute of Technology Delhi,\newline\indent New Delhi - 110016, India.}
\email{shravankumar.nageswaran@gmail.com}

\begin{document}


	\begin{abstract}
    
    Let $G$ a locally compact group and $(\Phi,\Psi)$ be a complementary pair of Young functions. Let $A_\Phi(G)$ be the Orlicz analogue of the classical Fig\`{a}-Talamanca Herz algebra $A_p(G).$ In this article, we establish a necessary and sufficient condition for the equality $\Lambda(A_\Phi(G)^{\ast\ast}) = A_\Phi(G)$ to hold, where $\Lambda(A_\Phi(G)^{\ast\ast})$ denotes the topological center of the double dual of $A_\Phi(G)$ when equipped with the first Arens product. Furthermore, we prove several results concerning the semi-simplicity of the Banach algebras $A_\Phi(G)^{\ast\ast}$ and $UCB_\Psi(\widehat{G})^\ast.$
    
	\end{abstract}
    
    
	\keywords{Amenability; Arens product; Banach algebras; Fig\`{a}-Talamanca Herz algebras; Orlicz spaces; Topological center.\\
		$\star$ Corresponding author. \\
		{\it Email address:} arvishdabra3@gmail.com (A. Dabra).}
	
	
	\maketitle

	
	\section{Introduction} 
    
    Let $\mathcal{A}$ be a Banach algebra and let $\mathcal{A^{\ast\ast}}$ denote its double dual. In 1951, Arens \cite{arens} introduced the first and second Arens products on $\mathcal{A}^{\ast\ast},$ proving that $\mathcal{A}^{\ast\ast}$ is a Banach algebra when equipped with either of these products. In this article, we refer to the first Arens product as the Arens product and denote it by $\square.$ From the definition, it follows that the map $\Tilde{\Gamma} \mapsto \Tilde{\Gamma} \, \square \, \Gamma$ is $w^\ast-$ continuous on $\mathcal{A}^{\ast\ast}$ for any $\Gamma \in \mathcal{A}^{\ast\ast}.$ In contrast, the map $\Tilde{\Gamma} \mapsto \Gamma \, \square \, \Tilde{\Gamma}$ is not necessarily $w^\ast-$ continuous. This asymmetry motivates the definition of the \textit{topological center} of $\mathcal{A}^{\ast\ast},$ given by
    $$\Lambda(\mathcal{A}^{\ast\ast}) := \{ \Gamma \in \mathcal{A}^{\ast\ast}: \Tilde{\Gamma} \mapsto \Gamma \, \square \, \Tilde{\Gamma} \,\, \text{is} \,\, w^\ast- \, \text{continuous on} \, \mathcal{A}^{\ast\ast} \}.$$ It is straightforward to verify that $\mathcal{A} \subseteq \Lambda(\mathcal{A}^{\ast\ast}).$

    Let $\mathcal{A}^{\ast}\mathcal{A}$ be the norm closure of the linear span of $\{T \cdot u: u \in \mathcal{A} \,\, \text{and} \,\, T \in \mathcal{A}^\ast\}$ in $\mathcal{A}^\ast.$ The dual space $(\mathcal{A}^{\ast}\mathcal{A})^\ast$ also becomes a Banach algebra when equipped with the multiplication induced by the Arens product on $\mathcal{A}^{\ast\ast}.$ Over the past three to four decades, the study of the topological center and related aspects of the algebras $\mathcal{A}^{\ast\ast}$ and $(\mathcal{A}^{\ast}\mathcal{A})^\ast$ has been the subject of considerable research interest. In 1987, Isik et al. \cite{INU} addressed the topological center problems for the group algebra $L^1(G)$ and established that \textit{the topological center of the algebra $L^1(G)^{\ast\ast}$ coincides with $L^1(G)$} for any compact group $G.$

    Earlier, in 1964, Eymard \cite{Eym} introduced and studied the Fourier algebra $A(G)$ and its dual $VN(G)$ for any locally compact group $G.$ If the underlying group $G$ is \textit{abelian} with $\widehat{G}$ as its dual group, then the Fourier algebra $A(G)$ is isometrically isomorphic to $L^1(\widehat{G}).$ Later, in 1974, Granirer \cite{Gra1974} defined the subspace of uniformly continuous functionals on $A(G),$ denoted by $UCB(\widehat{G}).$ When $G$ is abelian, this space coincides with the space of bounded uniformly continuous $\mathbb{C}-$valued functions on the dual group $\widehat{G}.$
    
    Cecchini and Zappa \cite{CZ} as well as Lau and Losert \cite{LL, LL2005}, investigated the \mbox{topological} center problems for the Banach algebra $\mathcal{A} = A(G)$ and the corresponding space $\mathcal{A}^\ast\mathcal{A} = UCB(\widehat{G}).$ Their results show that, for a large class of locally compact groups $G,$ including the abelian groups, discrete amenable groups, the motion group, the affine group $ax+b$ and the Heisenberg group, \textit{the topological center of the \mbox{algebra} $A(G)^{\ast\ast}$ is $A(G)$} and the topological center of $UCB(\widehat{G})^\ast$ is the Fourier-Stieltjes \mbox{algebra} $B(G).$ In 1996, Lau and \"Ulger \cite{LUl} extended this line of research further analyzing the topological center problems for both the group algebra $L^1(G)$ and the Fourier algebra $A(G),$ yielding significant results. In 2016, Losert \cite{VL} proved that if $G$ contains a free group $F_r \, (r \geq 2),$ then \textit{the topological center of $A(G)^{\ast\ast}$ is strictly larger than $A(G).$}

    For $1 < p < \infty$ and any locally compact group $G,$ Fig\`{a}-Talamanca \cite{figa} and Herz \cite{herz0} introduced the $L^p$-analogue of the Fourier algebra $A(G),$ known as the Fig\`{a}-Talamanca Herz algebra, denoted by $A_p(G).$ The algebra of $p-$pseudomeasures, denoted $PM_p(G),$ is defined as the $w^\ast-$ closure of $\lambda_p(M(G))$ in $\mathcal{B}(L^p(G)),$ where $\lambda_p$ is the representation of $M(G)$ on $L^p(G)$ given by left convolution. It is proved that the dual space of $A_p(G)$ is isometrically isomorphic to $PM_{p'}(G),$ where $p$ and $p'$ are dual H\"{o}lder exponents. In the special case $p = 2,$ the algebra $A_p(G)$ and its dual $PM_{p'}(G)$ coincide with the Fourier algebra $A(G)$ and the group von Neumann algebra $VN(G),$ respectively. 
     
    The algebra $A_p(G)$ has received significant attention in recent decades from numerous researchers, including Daws \cite{daws1}, Derighetti \cite{deriB, deriS}, Forrest \cite{FO1, FO2}, Granirer \cite{Gra1974, Gra0, Gra, GraS}, Lau and \"{U}lger \cite{lauUl} and Miao \cite{Miao1996, Miao2004, Miao2009AP, miao2009}, among others. The corresponding algebra of $p-$pseudomeasures $PM_p(G)$ has also been the subject of extensive study, notably in the work of Daws and Spronk \cite{daws2} as well as Gardella and Thiel \cite{GExpo, GH0, GH1, GH2}. In 2009, Miao \cite{miao2009} investigated the topological center problems for the algebra $\mathcal{A} = A_p(G)$ and the corresponding space $\mathcal{A}^\ast \mathcal{A} = UC_p(\widehat{G}),$ establishing a \textit{necessary and sufficient condition for the equality $\Lambda(A_p(G)^{\ast\ast}) = A_p(G)$ to hold}.

    It is well known that the Orlicz spaces are an important generalization of the \mbox{classical} Lebesgue spaces and have attracted considerable attention. Recently, Lal and Kumar \cite{RLSK1} and independently Aghababa and Akbarbaglu \cite{AA} introduced and studied the Orlicz $(L^\Phi)$-analogue of the Fig\`{a}-Talamanca Herz algebra $A_p(G),$ by replacing the $L^p$ spaces with $L^\Phi$ spaces, where $\Phi$ is a Young function \mbox{satisfying} the $\Delta_2-$condition. The resulting algebra, denoted $A_\Phi(G),$ is known as the \mbox{Orlicz} Fig\`{a}-Talamanca Herz algebra. If $(\Phi,\Psi)$ is a complementary pair of Young \mbox{functions}, then the dual space of $A_\Phi(G)$ is isometrically isomorphic to the space of $\Psi-$pseudomeasures, denoted by $PM_\Psi(G)$ \cite[Theorem 3.5]{RLSK1}. Both the algebra $A_\Phi(G)$ and its dual $PM_\Psi(G)$ have received a great deal of interest in recent years; see \cite{AA, AAM, AAR, AD, AD25, ArRaSh, RLSK1, RLSK3, RLSK2, RLSK4}.

    In this article, we investigate the topological center problems for the algebra $\mathcal{A} = A_\Phi(G)$ and the corresponding space $\mathcal{A}^\ast\mathcal{A} = UCB_\Psi(\widehat{G}),$ extending certain known results to a broader functional-analytic framework. The structure of the article is as follows: Section \ref{sec2} provides the necessary background, including basic definitions, notations and preliminary results. Section \ref{sec3} focuses on the existence of left and right identities in $A_\Phi(G)^{\ast\ast}$ and explore their relationship with the space $UCB_\Psi(\widehat{G}).$ In particular, Theorem \ref{rightidentity} characterizes the amenability of the group $G$ in terms of the existence of a right identity in $A_\Phi(G)^{\ast\ast},$ under the assumption that the Young function $\Phi$ satisfies the MA condition.
    
    In Section \ref{sec4}, we prove Theorem \ref{restriction}, which provides a sufficient condition under which the restriction of an element in $A_\Phi(G)^{\ast\ast}$ to $UCB_\Psi(\widehat{G})$ lies in $A_\Phi(G).$ This result plays a central role in establishing the main theorem of the article regarding the topological center of the double dual of $A_\Phi(G).$ Section \ref{sec5} is devoted to investigating the topological center problems for both the algebra $A_\Phi(G)$ and the corresponding space $UCB_\Psi(\widehat{G}).$ In Theorem \ref{topo}, we establish a necessary and sufficient condition for the equality $\Lambda(A_\Phi(G)^{\ast\ast}) = A_\Phi(G)$ to hold, assuming that $G$ is amenable. This result can be seen as an Orlicz analogue of Miao's theorem for $A_p(G)$ algebras. Moreover, an analogous result (Theorem \ref{lastthm}) is established for the topological center of the Banach algebra $UCB_\Psi(\widehat{G})^\ast.$


    
    \section{Preliminaries and Some Notations}\label{sec2}

    This section begins with a review of the fundamental terminology related to Orlicz spaces.

    An even convex function $ \Phi: \mathbb{R}\rightarrow [0,\infty]$ is termed as a \textit{Young function} if it satisfies the conditions $\Phi(0)= 0$ and $\lim\limits_{x \to \infty} \Phi(x)= + \infty.$ According to this definition, a Young function may attain the value $\infty$ at certain points and may therefore be discontinuous. However, in this context, we shall always assume that $\Phi$ is real-valued and continuous.
     
    For each Young function $\Phi,$ a corresponding convex function $\Psi,$ known as the \textit{complementary function} to $\Phi,$ is defined by
    $$\Psi(y):= \sup{\{x\, |y|-\Phi(x):x\geq 0\}}, \hspace{1cm} y\in\mathbb{R}.$$
    The function $\Psi$ thus defined is also a Young function. The pair $(\Phi,\Psi)$ (or equivalently $(\Psi,\Phi)$) is referred to as a \textit{complementary pair of Young functions}, or a complementary Young pair. As an example, for $1 < p < \infty,$ the function $\Phi(x) = |x|^p/p$ is a Young function, and its complementary function is $\Psi(y) = |y|^q/q,$ where $p$ and $q$ are dual H\"{o}lder exponents. 

    Let $G$ be a locally compact group equipped with a left Haar measure $dx.$ A Young function $\Phi$ is said to satisfy the $\Delta_{2}-$condition (globally), denoted $\Phi \in \Delta_{2} \, (\text{or} \,\, \Phi \in \Delta_2 \, \text{globally})$, if there exists a constant $k > 0$ and $x_{0} \geq 0$ such that $$\Phi(2x) \leq k \, \Phi(x)$$ for all $x \geq x_{0} \, (x_0 = 0)$ \cite[Definition 1, Pg. 22]{RR}. By \cite[Remark, Pg. 46]{RR}, $\Phi$ is said to be $\Delta_2-$\textit{regular} if the condition holds locally for $\Phi$ when the measure is finite and globally when the measure is infinite. A pair of complementary Young functions $(\Phi,\Psi)$ is said to satisfy the $\Delta_2-$condition if both $\Phi$ and $\Psi$ individually satisfy the $\Delta_2-$condition. It is straightforward to verify that the family of Young functions $\Phi_{\alpha,p}: x \mapsto \alpha \, |x|^p$ for $p \geq 1$ and $\alpha > 0,$ satisfies the $\Delta_2-$condition. Another example is the function $\Phi(x) = (1+|x|)\log(1+|x|)-|x|,$ which satisfies the $\Delta_2-$condition However, its complementary function $\Psi(y) = e^{|y|}-|y|-1$ does not satisfy the $\Delta_2-$condition. 
    
    A Young function $\Phi$ is said to satisfy the \textit{Milnes-Akimoni$\check{c}$} condition (abbreviated as the MA condition) if, for each $\epsilon > 0,$ there exists constants $\beta_\epsilon > 1$ and $x_0(\epsilon) \geq 0$ such that $$\Phi'((1+\epsilon)x) \geq \beta_\epsilon \, \Phi'(x), \hspace{1cm} x \geq x_0(\epsilon).$$ An example of a Young function satisfying the MA condition is $\Psi(y) = e^{|y|} - |y| - 1,$ for $y \in \mathbb{R}.$ However, its complementary function $\Phi(x) = (1+|x|)\log(1+|x|)-|x|,$ for $x \in \mathbb{R},$ does not satisfy the MA condition. It is easy to verify that the family of Young functions defined by $\Phi_\alpha(x) = |x|^\alpha (1 + |\log|x||),$ with $\alpha > 1,$ satisfies the MA condition. Moreover, the corresponding family $(\Phi_\alpha,\Psi_\alpha)$ of complementary pairs of Young functions also satisfies the $\Delta_2-$condition.
    
    For a given Young function $\Phi,$ the associated \textit{Orlicz space} is denoted by $L^{\Phi}(G)$ and is defined as
    $${L}^{\Phi}(G) := \left\{ f: G \rightarrow  \mathbb{C}:f \, \mbox{is measurable and}\int_G\Phi(\beta |f|)\ < \infty \text{ for some}~ \beta>0  \right\}.$$ The space $L^\Phi(G)$ becomes a Banach space when equipped with the Luxemburg norm (also known as the Gauge norm), defined by
    $$N_{\Phi}(f):= \inf \left\{k>0:\int_G\Phi\left(\frac{|f|}{k}\right) \leq1 \right\}.$$ Let $\Psi$ be the complementary function to $\Phi.$ Then, the Orlicz norm $\|\cdot\|_{\Phi}$ on $L^\Phi(G)$ is given by $$\|f\|_{\Phi} := \sup \left\{\int_{G}|fg| :g\in L^\Psi(G) \, \, \text{and} \, \, \int_{G}\Psi(|g|) \leq1 \right\}.$$ In fact, these two norms are equivalent and for every $f \in L^\Phi(G),$ the following holds: $$N_\Phi(f) \leq \|f\|_\Phi \leq 2 \, N_\Phi(f).$$ 

    If a Young function $\Phi$ satisfies the $\Delta_{2}-$condition, then the space $\mathcal{C}_c(G)$ of all continuous functions on $G$ with compact support is dense in $L^\Phi(G).$ Further, if the complementary function $\Psi$ is continuous and satisfies $\Psi(x) = 0$ if and only if $x = 0,$ then the dual space of $(L^\Phi(G),N_\Phi(\cdot))$ is isometrically isomorphic to $(L^\Psi(G),\|\cdot\|_\Psi)$ \cite[Corollary 9, Pg. 111]{RR}. In particular, if the complementary Young pair $(\Phi,\Psi)$ satisfies the $\Delta_2-$condition, then the Banach space $L^\Phi(G)$ is reflexive \cite[Theorem 10, Pg. 112]{RR}.
    
    For more detailed information on Orlicz spaces, we refer the readers to \cite{RR}.
    \newline

    Let $\Phi$ be a Young function satisfying the $\Delta_2-$condition. For any function $g:G\to\mathbb{C}$, we define $\check{g}$ by $\check{g}(x):= g(x^{-1})$ for all $x \in G.$ The \textit{Orlicz Fig\`a-Talamanca Herz algebra} $A_\Phi(G)$ is defined as the set of all continuous functions $u \in \mathcal{C}_0(G)$ that can be expressed as $$u=\sum_{n \in \mathbb{N}} f_n \ast \check{g_n},$$ where $f_n\in L^\Phi(G)$ and $g_n\in L^\Psi(G),$ such that $$\sum_{n \in \mathbb{N}} N_\Phi(f_n)\|g_n\|_\Psi<\infty.$$ For any $u \in A_\Phi(G),$ the norm is defined by $$\|u\|_{A_\Phi} := \inf\left\{\sum\limits_{n \in \mathbb{N}} N_\Phi(f_n) \|g_n\|_\Psi: u = \sum\limits_{n \in \mathbb{N}} f_n \ast \check{g_n}\right\}.$$ Equipped with this norm, along with pointwise addition and multiplication, the space $A_\Phi(G)$  becomes a commutative Banach algebra. Moreover, it is a regular, tauberian, semi-simple Banach algebra with its Gelfand spectrum homeomorphic to $G.$

    Recall from \cite[Lemma 4.1]{RLSK1} that for any open subgroup $H$ of $G,$ the Banach algebra $A_\Phi(H)$ is isometrically isomorphic to the subalgebra of $A_\Phi(G)$ consisting of functions that vanish outside $H.$ We denote by $i_H$ the \textit{inclusion map} from $A_\Phi(H)$ to $A_\Phi(G)$  and  by $m_H$ the \textit{restriction map} from $A_\Phi(G)$ to $A_\Phi(H).$ It follows from \cite[Lemma 4.5]{RLSK1} that the map $m_H$ is a contractive homomorphism of Banach algebras. Further, it is easy to verify that $A_\Phi(H)$ is a closed ideal of $A_\Phi(G)$ using the aforementioned identification. Moreover, the map $m_H$ is a bounded projection which is also a multiplier, i.e., it satisfies $$m_H(u\,v) = u \, m_H(v) = m_H(u) \, v$$ for all $u,v \in A_\Phi(G)$ and $m_H(u) = u$ for all $u \in A_\Phi(H).$

    Let $B_\Phi(G)$ denote the multiplier algebra of $A_\Phi(G),$ defined as the set $$\{u \in \mathcal{C}(G): u \, v \in A_\Phi(G) \,\, \text{for all} \,\, v \in A_\Phi(G)\}.$$ Each element $u \in B_\Phi(G)$ induces a bounded linear operator $M_u: A_\Phi(G) \to A_\Phi(G),$ given by $M_u(v) := u \, v,$ for all $v \in A_\Phi(G).$ Equipped with the operator norm and pointwise operations (addition and multiplication), the space $B_\Phi(G)$ is a commutative Banach algebra. Moreover, it is clear that $A_\Phi(G) \subseteq B_\Phi(G).$

    Let $M(G)$ denote the set of all bounded complex Radon measures on $G.$ For $\mu \in M(G)$ and $f \in L^\Psi(G),$ we define the operator $T_\mu: L^\Psi(G) \to L^\Psi(G)$ by left convolution: $$T_\mu(f) := \mu \ast f.$$ If $\mathcal{B}(L^\Psi(G))$ denote the Banach space of all bounded linear operators on $L^\Psi(G),$ equipped with the operator norm. It is straightforward to verify that $T_\mu \in \mathcal{B}(L^\Psi(G))$ for each $\mu \in M(G).$ Let $PM_\Psi(G)$ be the closure of the set $\{T_\mu: \mu \in M(G)\}$ in $\mathcal{B}(L^\Psi(G))$ with respect to the ultra-weak topology (the $w^\ast$-topology). It is proved that the dual space of $A_\Phi(G)$ is isometrically isomorphic to $PM_\Psi(G)$ \cite{RLSK1}.
    
    Let $PF_\Psi(G)$ denote the norm closure in $\mathcal{B}(L^\Psi(G))$ of the set $\{T_{\mu_f}: f \in L^1(G)\},$ where $\mu_f$ is the measure associated with $f \in L^1(G).$ The dual of $PF_\Psi(G)$ is denoted by $W_\Phi(G),$ which is a commutative Banach algebra containing $A_\Phi(G)$ \cite{RLSK4}. For each $\mu = \delta_x \, \, (\text{where} \,\, x \in G),$ the associated operator $T_{\delta_x}$ is denoted by $\lambda_\Psi(x).$ The map $\lambda_\Psi$ thus defines the left regular representation of $G$ on $L^\Psi(G).$

    A linear functional $\Gamma \in PM_\Psi(G)^\ast$ is called a \textit{topologically invariant mean} if it satisfies the following conditions: $$\|\Gamma\| = 1 = \langle \Gamma,\lambda_\Psi(e)\rangle \,\, \text{and} \,\, \langle \Gamma,u \cdot T \rangle = u(e) \, \langle \Gamma,T \rangle$$ for all $u \in A_\Phi(G)$ and $T \in PM_\Psi(G).$ We denote the set of all topologically invariant means on $PM_\Psi(G)$ by $TIM_\Psi(\widehat{G}).$

    For additional details on the Banach algebra $A_\Phi(G)$ and the related spaces, we refer the readers to the series of papers \cite{AA, AAR, AD, AD25, ArRaSh, RLSK1, RLSK3, RLSK2, RLSK4}.
    \newline
    
    Let $\mathcal{A}$ be a Banach algebra and let $\mathcal{A}^{\ast\ast}$ denote its double dual. It is well known that $\mathcal{A}^{\ast\ast}$ can be equipped with two natural extensions of the multiplication on $\mathcal{A},$ both of which turn $\mathcal{A}^{\ast\ast}$ into a Banach algebra. These extensions, introduced by Arens \cite{arens}, are commonly referred to as the first and second Arens product. In this article, we focus on \textit{the first Arens product}, denoted by $\square,$ which is defined as follows:
	\begin{enumerate}[(i)]
		\item $\langle u \cdot T,v\rangle = \langle T,u \, v\rangle \,\, \text{for every} \,\, u,v \in \mathcal{A} \,\, \text{and} \,\, T \in \mathcal{A}^\ast.$
		\item $\langle \Gamma \odot T,u\rangle = \langle \Gamma,u \cdot T\rangle \,\, \text{for every} \,\, u \in \mathcal{A}, T \in \mathcal{A}^\ast \,\, \text{and} \,\, \Gamma \in \mathcal{A}^{\ast\ast}.$
		\item $\langle \Tilde{\Gamma} \, \square \, \Gamma,T\rangle = \langle \Tilde{\Gamma},\Gamma \odot T\rangle \,\, \text{for every} \,\, T \in \mathcal{A}^\ast \,\, \text{and} \,\, \Tilde{\Gamma},\Gamma \in \mathcal{A}^{\ast\ast}.$
	\end{enumerate}

    Recall that a Banach algebra is said to be \textit{weakly sequentially complete} (WSC) if every weakly Cauchy sequence converges with respect to the weak topology.\\

    A locally compact group $G$ is said to be \textit{amenable} if there exists a left-invariant mean on $L^\infty(G),$ i.e., a positive linear functional $\Lambda$ on $L^\infty(G)$ of norm one satisfying $$\langle \Lambda,L_x f \rangle = \langle \Lambda,f \rangle \,\, \text{for all} \,\, f \in L^\infty(G) \,\, \text{and} \,\, x \in G,$$ where $L_xf(y) := f(x^{-1}y)$ denotes the left translation of $f.$ Classical examples of amenable groups include abelian groups, compact groups and solvable groups. In contrast, the free group on two generators is a well known example of a non-amenable group.

    Throughout this article, $G$ denotes a locally compact group with a fixed Haar measure $dx$ and assume that the complementary pair of Young functions $(\Phi,\Psi)$ satisfies the $\Delta_2-$condition.



    \section{The Double Dual of $A_\Phi(G)$ and the Space $UCB_\Psi(\widehat{G})$}\label{sec3}
    
    In this section, we present certain results related to the left and right identities in the double dual of $A_\Phi(G)$ when considered as a Banach algebra with the first Arens product. Theorem \ref{rightidentity} establishes a necessary and sufficient condition for the existence of a right identity in $A_\Phi(G)^{\ast\ast}.$ Further, we explore their relationship with the space $UCB_\Psi(\widehat{G}).$ Moreover, the semi-simplicity of the Banach algebras $A_\Phi(G)^{\ast\ast}$ and $UCB_\Psi(\widehat{G})^\ast$ is studied in Theorem \ref{dualsemi} and Theorem \ref{ucbsemi}, respectively.
    
    \begin{rem}
    We would like to mention that, to the best of our knowledge, the results presented in this section are, in fact, new for the Fig\`{a}-Talamanca Herz algebra $A_p(G),$ i.e., when the corresponding Young function $\Phi(x)$ is given by $|x|^p/p,$ for $1 < p < \infty.$
    \end{rem}

    Let us begin with the following propositions related to the existence of left identity in the Banach algebra $A_\Phi(G)^{\ast\ast}.$

    \begin{prop}
    The algebra $A_\Phi(G)^{\ast\ast}$ has a left identity if $G$ is compact.
    \end{prop}

    \begin{proof}
    Since $G$ is compact, the algebra $A_\Phi(G)$ contains the unit function $\textbf{1}$ which also serves as the identity. Then, it is easy to verify by elementary calculations and using the definition of the Arens product $\square$ that $J(\textbf{1})$ is a left identity for $A_\Phi(G)^{\ast\ast},$ where $J$ is the canonical map from $A_\Phi(G)$ to $A_\Phi(G)^{\ast\ast}.$ 
    \end{proof}

    \begin{prop}\label{leftidentity}
    If $A_\Phi(G)^{\ast\ast}$ has a left identity, then $UCB_\Psi(\widehat{G}) = PM_\Psi(G).$
    \end{prop}

    \begin{proof}
    If possible, assume that $UCB_\Psi(\widehat{G})$ is a proper closed subspace of $PM_\Psi(G).$ Then, by Hahn-Banach theorem, there exists a non-zero $\Gamma \in PM_\Psi(G)^\ast = A_\Phi(G)^{\ast\ast}$ such that $\Gamma(UCB_\Psi(\widehat{G})) = \{0\}.$ This implies that $\Gamma \odot T = 0$ for all $T \in PM_\Psi(G)$ and hence, $\Tilde{\Gamma} \, \square \, \Gamma = 0$ for all $\Tilde{\Gamma} \in A_\Phi(G)^{\ast\ast}.$ Since $A_\Phi(G)^{\ast\ast}$ has a left identity, say $\Gamma_e,$ we have, $$\Gamma = \Gamma_e \, \square \, \Gamma = 0,$$ which is a contradiction. Thus, $UCB_\Psi(\widehat{G}) = PM_\Psi(G).$
    \end{proof}

    Here is the main theorem of this section, which characterizes the amenability of the group $G$ in terms of the existence of a right identity in the algebra $A_\Phi(G)^{\ast\ast}.$ This result generalizes Lau's result \cite[Proposition 3.2 (a)]{Lau} on $VN(G)^\ast$ for any locally compact group $G.$ Furthermore, the analogous result for hypergroups is given by Esmailvandi et al. \cite[Proposition 4.1]{RNS}.

    \begin{thm}\label{rightidentity}
    Let $\Phi$ satisfies the MA condition. Then the Banach algebra $A_\Phi(G)^{\ast\ast}$ has a right identity if and only if $G$ is amenable.
    \end{thm}
    
    \begin{proof}
    Assume that $G$ is amenable. Then, by \cite[Theorem 3.1]{RLSK2}, the algebra $A_\Phi(G)$ possesses a bounded approximate identity. Hence, by \cite[Proposition III.28.7]{BD}, there exists a right identity in the Banach algebra $A_\Phi(G)^{\ast\ast}.$

    Conversely, let $\Gamma_e$ be a right identity in $A_\Phi(G)^{\ast\ast} = PM_\Psi(G)^\ast.$ Now, by Goldstine's theorem (\cite[Theorem A.3.29]{dales}), there exists a net $(u_\alpha)$ in $A_\Phi(G)$ such that $$\|u_\alpha\| \leq \|\Gamma_e\| \,\, \text{and} \,\, \Gamma_e = w^\ast-\lim_\alpha J(u_\alpha),$$
    where $J$ is the canonical map from $A_\Phi(G)$ to $A_\Phi(G)^{\ast\ast}.$
        
    Then, for every $u \in A_\Phi(G)$ and $T \in PM_\Psi(G),$ 
    $$\langle T,u \rangle = \langle J(u),T \rangle = \langle J(u) \, \square \, \Gamma_e, T \rangle = \lim_\alpha \, \langle T,u_\alpha u \rangle.$$
    This implies that $(u_\alpha)$ is a bounded weak approximate identity for $A_\Phi(G).$ Thus, by \cite[Proposition I.11.4]{BD}, the algebra $A_\Phi(G)$ has a bounded approximate identity and hence, by \cite[Theorem 3.1]{RLSK2}, $G$ is amenable. 
    \end{proof}

    In the subsequent results, we investigate the semi-simplicity of the algebras $A_\Phi(G)^{\ast\ast}$ and $UCB_\Psi(\widehat{G})^\ast,$ establishing certain
    necessary and sufficient conditions.

    \begin{prop}
    If the group $G$ is finite, then the algebra $A_\Phi(G)^{\ast\ast}$ is semi-simple.
    \end{prop}

    \begin{proof}
    The result follows from \cite[Theorem 3.8]{AD} and the fact that the algebra $A_\Phi(G)$ is semi-simple.
    \end{proof}

    The following theorem demonstrates that the semi-simplicity of the algebra $A_\Phi(G)^{\ast\ast}$ is sufficient for a group $G$ to be discrete, provided that $G$ is a second countable group for which the algebra $A_\Phi(G)$ is WSC. The precise statement of the result is as follows.

    \begin{thm}\label{dualsemi}
    Let $G$ be a second countable group such that the algebra $A_\Phi(G)$ is weakly sequentially complete. If $A_\Phi(G)^{\ast\ast}$ is semi-simple, then $G$ is discrete.
    \end{thm}

    \begin{proof}
    Let 
    $$I := \left\{\begin{array}{c}
    \Gamma \in A_\Phi(G)^{\ast\ast}: \langle \Gamma,\lambda_\Psi(e) \rangle = 0 \,\, \text{and} \,\, \langle \Gamma,u \cdot T \rangle = u(e) \, \langle \Gamma,T \rangle \\ \forall \, u \in A_\Phi(G) \,\, \text{and} \,\, T \in PM_\Psi(G)
    \end{array} \right\}.$$
    We claim that $I$ is an ideal in $A_\Phi(G)^{\ast\ast}.$

    Observe that, for any $\Tilde{\Gamma} \in A_\Phi(G)^{\ast\ast}$ and $\Gamma \in I,$
    $$\langle \Tilde{\Gamma} \, \square \, \Gamma,\lambda_\Psi(e) \rangle = \langle \Tilde{\Gamma},\Gamma \odot \lambda_\Psi(e) \rangle = \langle \Tilde{\Gamma},\langle \Gamma,\lambda_\Psi(e) \rangle \, \lambda_\Psi(e) \rangle = \langle \Gamma,\lambda_\Psi(e) \rangle \, \langle \Tilde{\Gamma},\lambda_\Psi(e) \rangle = 0.$$

    Since, the condition $\langle \Gamma,u \cdot T \rangle = u(e) \, \langle \Gamma,T \rangle$ is equivalent to $\Gamma \, \square \, u = u(e) \, \Gamma$ and the Arens product $\square$ is associative, it follows that
    $$(\Tilde{\Gamma} \, \square \, \Gamma)\, \square \, u = \Tilde{\Gamma} \, \square \, (\Gamma \, \square \, u) = \Tilde{\Gamma} \, \square \, (u(e) \, \Gamma) = u(e) \, (\Tilde{\Gamma} \, \square \, \Gamma).$$
    This implies $\Tilde{\Gamma} \, \square \,  \Gamma \in I$ for all $\Gamma \in I$ and $\Tilde{\Gamma} \in A_\Phi(G)^{\ast\ast}.$

    Further, by using the fact that the algebra $A_\Phi(G)$ is contained in the center of $A_\Phi(G)^{\ast\ast} = PM_\Psi(G)^\ast$ \cite[Theorem 3.2]{ArRaSh}, one can similarly prove that $\Gamma \, \square \, \Tilde{\Gamma} \in I.$ Thus, $I$ is an ideal in $A_\Phi(G)^{\ast\ast}.$

    Moreover, for any $\Tilde{\Gamma}, \Gamma \in I$ and $T \in PM_\Psi(G),$
    $$\langle \Tilde{\Gamma} \, \square \, \Gamma, T \rangle = \langle \Tilde{\Gamma}, \Gamma \odot T \rangle = \langle \Tilde{\Gamma}, \langle \Gamma,T \rangle \, \lambda_\Psi(e) \rangle = \langle \Gamma,T \rangle \, \langle \Tilde{\Gamma},\lambda_\Psi(e) \rangle = 0.$$
    This implies that $I^2 = \{0\},$ i.e., $I$ is nil and thus, by \cite[Proposition 1.5.6 (ii)]{dales}, 
    $$I \subseteq rad(A_\Phi(G)^{\ast\ast}).$$

    Now, if possible, assume that $\Gamma_1$ and $\Gamma_2$ are two distinct elements of $TIM_\Psi(\widehat{G}).$ Then, by definition, $0 \neq \Gamma_1 - \Gamma_2 \in I,$ which contradicts the fact that the Banach algebra $A_\Phi(G)^{\ast\ast}$ is semi-simple. This implies that $TIM_\Psi(\widehat{G})$ is singleton and hence, by \cite[Corollary 7.6]{RLSK1}, $G$ is discrete.
    \end{proof}

    \begin{prop}
    If $A_\Phi(G)^{\ast\ast}$ is semi-simple, then $UCB_\Psi(\widehat{G}) = PM_\Psi(G).$
    \end{prop}

    \begin{proof}
    Let us suppose to the contrary that $UCB_\Psi(\widehat{G})$ is a proper closed subspace of $PM_\Psi(G).$ Let $K$ be the annihilator of $UCB_\Psi(\widehat{G})$ in $A_\Phi(G)^{\ast\ast},$ i.e.,
    $$K = \{\Gamma \in A_\Phi(G)^{\ast\ast}: \Gamma(UCB_\Psi(\widehat{G})) = \{0\} \}.$$
    Then, by Hahn-Banach theorem, $K \neq \{0\}.$ Further, it is clear from the proof of Proposition \ref{leftidentity} that $K$ is a left ideal in $A_\Phi(G)^{\ast\ast}$ and $K^2 = \{0\}.$ In particular, we have, $K \subseteq rad(A_\Phi(G)^{\ast\ast})$ and hence, this contradicts the semi-simplicity of the algebra $A_\Phi(G)^{\ast\ast}.$ Thus, $UCB_\Psi(\widehat{G}) = PM_\Psi(G).$
    \end{proof}

    \begin{rem}\label{semi}
    Note that the Banach algebras $W_\Phi(G)$ and $B_\Phi(G)$ are semi-simple. This follows directly from the relation $A_\Phi(G) \subseteq W_\Phi(G) \subseteq B_\Phi(G)$ \cite[Corollary 3.2]{RLSK4} and the fact that the set $\{ \lambda_\Psi(x) : x \in G\}$ is contained in the Gelfand spectrum of all of these algebras.
    \end{rem}

    The following theorem addresses the semi-simplicity of the algebra $UCB_\Psi(\widehat{G})^\ast$ and is analogous to Theorem \ref{dualsemi}.
    
    \begin{thm}\label{ucbsemi}
    The Banach algebra $UCB_\Psi(\widehat{G})^\ast$ is semi-simple if $G$ is discrete. The converse holds when $G$ is a second countable group for which $A_\Phi(G)$ is weakly sequentially complete.
    \end{thm}

    \begin{proof}
    Since $G$ is discrete, by \cite[Corollary 4.13 (ii)]{RLSK3}, $UCB_\Psi(\widehat{G}) = PF_\Psi(G)$ and thus, $UCB_\Psi(\widehat{G})^\ast = PF_\Psi(G)^\ast = W_\Phi(G)$ \cite{RLSK4}. The desired result follows directly from Remark \ref{semi}.

    Conversely, let $G$ be a second countable group such that the algebra $A_\Phi(G)$ is weakly sequentially complete and assume that $UCB_\Psi(\widehat{G})^\ast$ is semi-simple. Consider
    $$I := \left\{\begin{array}{c}
    \Gamma \in UCB_\Psi(\widehat{G})^\ast: \langle \Gamma,\lambda_\Psi(e) \rangle = 0 \,\, \text{and} \,\, \langle \Gamma,u \cdot T \rangle = u(e) \, \langle \Gamma,T \rangle \\ \forall \, u \in A_\Phi(G) \,\, \text{and} \,\, T \in UCB_\Psi(\widehat{G})
    \end{array} \right\}.$$
    Now, by repeating the arguments as in Theorem \ref{dualsemi}, it follows that $I$ is an ideal in $UCB_\Psi(\widehat{G})^\ast$ and $I \subseteq rad(UCB_\Psi(\widehat{G})^\ast).$ Since, by \cite[Corollary 6.2]{RLSK1}, $TIM_\Psi(\widehat{G}) \neq \emptyset,$ let $\Gamma_r$ be the restriction of $\Gamma$ to $UCB_\Psi(\widehat{G}),$ for some non-zero $\Gamma \in TIM_\Psi(\widehat{G}).$ Observe that $\Gamma_r$ is also non-zero and belongs to $UCB_\Psi(\widehat{G})^\ast.$ Otherwise, if $\Gamma_r = 0,$ then $$\langle \Gamma, u \cdot T \rangle = \langle \Gamma_r, u \cdot T \rangle = 0 \,\, \text{for all} \,\, u \in A_\Phi(G) \,\, \text{and} \,\, T \in PM_\Psi(G).$$
    As $\Gamma \in TIM_\Psi(\widehat{G}),$ 
    $$0 = \langle \Gamma,u \cdot T \rangle = u(e) \, \langle \Gamma,T \rangle \,\, \text{for all} \,\, u \in A_\Phi(G) \,\, \text{and} \,\, T \in PM_\Psi(G).$$
    Now, by choosing a suitable $u \in A_\Phi(G)$ such that $u(e) = 1,$ we have,
    $$\langle \Gamma,T \rangle = 0 \,\, \text{for all} \,\, T \in PM_\Psi(G).$$
    This contradicts the fact that $\Gamma$ is non-zero and hence, $\Gamma_r \neq 0.$ Thus, by repeating the steps as in Theorem \ref{dualsemi}, it follows that $TIM_\Psi(\widehat{G})$ is a singleton set and hence, by \cite[Corollary 7.6]{RLSK1}, $G$ is discrete.
    \end{proof}
    

    
    \section{Restriction Theorem}\label{sec4}

    The aim of this section is to prove the restriction theorem (Theorem \ref{restriction}), which plays a pivotal role in establishing the main result of this article (Theorem \ref{topo}), concerning the topological center of the double dual of $A_\Phi(G).$ The ideas behind the proofs are motivated by \cite{miao2009}.

    We begin with the following proposition.
    
    \begin{prop}
    Let $H$ be an open subgroup of an amenable group $G.$ Then every right unit of $A_\Phi(H)^{\ast\ast}$ in $A_\Phi(G)^{\ast\ast}$ admits an extension to a right unit of $A_\Phi(G)^{\ast\ast}.$
    \end{prop}

    \begin{proof}
    Since $G$ is amenable, by \cite[Theorem 3.1]{RLSK2}, the algebra $A_\Phi(G)$ has a bounded approximate identity. Further, as the restriction map $m_H: A_\Phi(G) \to A_\Phi(H)$ is a bounded projection which is also a multiplier, the result is a direct consequence of \cite[Theorem 2.3]{miao2009}.
    \end{proof}

    Recall from \cite[Pg. 101]{herz} that for an algebra of functions $\mathcal{A},$ the support of a linear functional $T \in \mathcal{A}^\ast$ as a subset of $G$ is characterized as follows: $x \notin supp(T)$ if and only if there exists a neighborhood $V_x$ of $x$ such that $\langle T,v \rangle=0$ for all $v \in \mathcal{A}$ with $supp(v) \subseteq V_x.$

    \begin{rem}\label{supp}
    It is easy to verify that if $supp(v) \cap supp(T) = \emptyset$ for $v \in A_\Phi(G)$ and $T \in PM_\Psi(G),$ then the function $v \cdot T = 0.$
    \end{rem}

    \begin{lem}\label{open}
    Let $H$ be an open subgroup of $G$ and $m_H: A_\Phi(G) \to A_\Phi(H)$ be the restriction map. Then $m_H^\ast(T) = T$ for all $T \in PM_\Psi(G)$ with $supp(T) \subseteq H.$
    \end{lem}

    \begin{proof}
    Let $u,v \in A_\Phi(G) \cap C_c(G)$ be such that $supp(u) \cap (G \setminus H) \neq \emptyset$ and $v(x) = 1$ for all $x \in supp(u) \cap (G\setminus H)$ with $supp(v) \subseteq G \setminus H.$ Then, for any $T \in PM_\Psi(G)$ with $supp(T) \subseteq H,$ we have,
    $$\langle T,u \rangle  = \langle T,u|_H \rangle + \langle T,u \, v \rangle = \langle T,m_H(u) \rangle + \langle v \cdot T,u \rangle.$$
    Now, by Remark \ref{supp}, it follows that
    $$\langle T,u \rangle = \langle m_H^\ast(T),u \rangle \,\, \text{for all} \,\, u \in A_\Phi(G) \cap C_c(G).$$
    Thus, by the density of $A_\Phi(G) \cap C_c(G)$ in $A_\Phi(G),$ it follows that $m_H^\ast(T) = T$ for all $T \in PM_\Psi(G)$ with $supp(T) \subseteq H.$
    \end{proof}

    The next lemma is a crucial component in proving Theorem \ref{restriction} and establishes the existence of a sequence of compact subsets of $G$ with certain significant properties.
    
    \begin{lem}\label{existence}
    Let $\Gamma \in A_\Phi(G)^{\ast\ast}$ and assume that for every open $\sigma-$compact subgroup $H$ of $G,$ $m_H^{\ast\ast}(\Gamma)$ is in $A_\Phi(H).$ Then, for each $n \in \mathbb{N},$ there exists a compact subset $K_n$ of $G$ such that $|\langle \Gamma,T \rangle| < 1/n$ for all $T \in UCB_\Psi(\widehat{G})$ with compact support, $\|T\| \leq 1$ and $supp(T) \subseteq G \setminus K_n.$
    \end{lem}

    \begin{proof}
    Let us suppose to the contrary that there exists $n_0 \in \mathbb{N}$ and a function $T_1 \in UCB_\Psi(\widehat{G})$ with compact support such that $\|T_1\| \leq 1$ and $|\langle \Gamma,T_1\rangle| \geq \epsilon,$ where $\epsilon = 1/n_0 > 0.$ Now, by reiterating the construction outlined in \cite[Lemma 3.1]{miao2009}, we obtain a sequence of compactly supported functions $\{T_n\}$ in $UCB_\Psi(\widehat{G})$ and a sequence of relatively compact symmetric neighborhoods $\{V_n\}$ of $e$ in $G$ such that, for each $n \in \mathbb{N},$
    \begin{enumerate}
            \item $supp(T_{n+1}) \subseteq G \setminus \overline{V_n}$ and $V_n^2 \subseteq V_{n+1};$
            \item $\|T_n\| \leq 1$ and $supp(T_n) \subseteq V_n;$
            \item $|\langle \Gamma,T_n\rangle| \geq \epsilon.$
    \end{enumerate}
    Let $H = \bigcup_n V_n.$ Then, $H$ is an open $\sigma-$compact subgroup of $G$ and by hypothesis, $m_H^{\ast\ast}(\Gamma)$ is in $A_\Phi(H).$ This implies that $m_H^{\ast\ast}(\Gamma) = J(v_\Gamma),$ for some $v_\Gamma \in A_\Phi(H),$ where $J$ is the canonical map from $A_\Phi(H)$ to $A_\Phi(H)^{\ast\ast}.$ Now, we claim that there exists a compact subset $K$ of $H$ such that $|\langle m_H^{\ast\ast}(\Gamma),T\rangle| \leq \epsilon/2$ for any $T \in PM_\Psi(G)$ with $\|T\| \leq 1$ and $supp(T) \subseteq G \setminus K.$ Since $A_\Phi(H) \cap C_c(H)$ is dense in $A_\Phi(H),$ there exists $u \in A_\Phi(H) \cap C_c(H)$ such that 
    $$\|v_\Gamma - u\| \leq \epsilon/2.$$
    Let $U$ be a relatively compact neighborhood of $supp(u)$ in $H.$ Then, the required compact set $K$ is $\overline{U}$ and observe that, for any $T \in PM_\Psi(G)$ with $\|T\| \leq 1$ and $supp(T) \subseteq G \setminus K,$ we have,
    \begin{align*}
            |\langle m_H^{\ast\ast}(\Gamma),T \rangle| = |\langle J(v_\Gamma),T \rangle| &= |\langle T,v_\Gamma \rangle| \\ &\leq |\langle T,v_\Gamma - u \rangle| + |\langle T,u\rangle| \\ &\leq \|v_\Gamma - u\| \leq \epsilon/2.
    \end{align*}

    In the above computation, since $supp(u) \cap supp(T) = \emptyset,$ the term $|\langle T,u \rangle|$ equals zero \cite[Pg. 101]{herz}. Now,
        $$K \subseteq H = \bigcup_n V_n,$$
        i.e., $\{V_n\}$ is an open cover of $K.$ Since, $K$ is compact and the sequence $\{V_n\}$ is increasing, there exists $k \in \mathbb{N}$ such that $K \subseteq V_k.$ Then, $$supp(T_{k+1}) \subseteq G \setminus \overline{V_k} \subseteq G \setminus K$$ and thus, by the above claim, it follows that $|\langle m_H^{\ast\ast}(\Gamma),T_{k+1}\rangle| \leq \epsilon/2.$ However, by Lemma \ref{open}, we have,
        $$|\langle m_H^{\ast\ast}(\Gamma),T_{k+1} \rangle| = |\langle\Gamma,m_H^\ast(T_{k+1})\rangle| = |\langle \Gamma,T_{k+1}\rangle| \geq \epsilon,$$
        which is a contradiction. Hence, the result follows.
    \end{proof}

    Here is the main theorem of this section, which provides a sufficient condition for any $\Gamma \in A_\Phi(G)^{\ast\ast},$ under which the restriction of $\Gamma$ to $UCB_\Psi(\widehat{G})$ belongs to $A_\Phi(G).$ The analogous result for the Fourier algebra $A(H),$ where $H$ is a hypergroup and for Fig\`a-Talamanca Herz algebra $A_p(G),$ where $G$ is a locally compact group, are given by Esmailvandi et al. \cite[Proposition 5.5]{RNS} and Miao \cite[Lemma 3.1]{miao2009}, respectively. The following theorem is the Orlicz analogue of Miao's result.
    
    \begin{thm}\label{restriction}
        Let $\Gamma \in A_\Phi(G)^{\ast\ast}.$ If for every open $\sigma-$compact subgroup $H$ of $G,$ $m_H^{\ast\ast}(\Gamma)$ is in $A_\Phi(H),$ then the restriction of $\Gamma$ to $UCB_\Psi(\widehat{G})$ belongs to $A_\Phi(G).$
    \end{thm}

    \begin{proof}
        Let $H$ be an open $\sigma-$compact subgroup of $G$ containing all compact subsets $K_n$ of $G$ from Lemma \ref{existence}. Then, for any $T \in UCB_\Psi(\widehat{G})$ with compact support, it is easy to verify that the function $(T - m_H^\ast(T))$ also belongs to $UCB_\Psi(\widehat{G})$ and has compact support contained in $G \setminus H.$ Since, for each $n \in \mathbb{N},$ $(G \setminus H) \subseteq (G \setminus K_n),$ by Lemma \ref{existence}, it follows that $\langle \Gamma,T - m_H^\ast(T)\rangle = 0.$ Thus,
        $$\langle \Gamma, T \rangle = \langle \Gamma, T - m_H^\ast(T)\rangle + \langle \Gamma,m_H^\ast(T)\rangle = \langle m_H^{\ast\ast}(\Gamma),T\rangle,$$
        for all $T \in UCB_\Psi(\widehat{G})$ with compact support and hence, by \cite[Corollary 4.7]{RLSK3}, the restriction of $\Gamma$ to $UCB_\Psi(\widehat{G})$ is $m_H^{\ast\ast}(\Gamma).$ The desired result is a direct consequence of the hypothesis that $m_H^{\ast\ast}(\Gamma)$ is in $A_\Phi(H) \subseteq A_\Phi(G).$
    \end{proof}



    \section{Topological center}\label{sec5}
    
    In this section, we study the topological center of the double dual of $A_\Phi(G),$ when considered as a Banach algebra with the first Arens product. The primary aim is to prove Theorem \ref{topo}, which establishes a necessary and sufficient condition for the equality $\Lambda(A_\Phi(G)^{\ast\ast}) = A_\Phi(G)$ to hold, assuming that $G$ is amenable. The ideas and results are inspired by Miao's work on $A_p(G)$ algebras \cite{miao2009}.
    
    Let us begin with a straightforward yet significant lemma.

    \begin{lem}\label{ihmh}
        Let $H$ be an open subgroup of $G$ and $i_H: A_\Phi(H) \to A_\Phi(G)$ be the inclusion map. Then the following holds
        $$i_H^{\ast\ast}(\Lambda(A_\Phi(H)^{\ast\ast})) \subseteq \Lambda(A_\Phi(G)^{\ast\ast}) \,\, \text{and} \,\,\, m_H^{\ast\ast}(\Lambda(A_\Phi(G)^{\ast\ast})) \subseteq \Lambda(A_\Phi(H)^{\ast\ast}).$$
    \end{lem}

    \begin{proof}
        Since the restriction map $m_H: A_\Phi(G) \to A_\Phi(H)$ is a bounded projection which is also a multiplier, the result is a direct consequence of \cite[Theorem 4.2]{miao2009}.
    \end{proof}

    Here is the main result of this article concerning the topological center of the double dual of $A_\Phi(G).$ This theorem serves as an Orlicz analogue of Miao's result on $A_p(G)$ algebras \cite[Theorem 4.4]{miao2009}.
    
    \begin{thm}\label{topo}
        Let $G$ be an amenable group and suppose that $\Phi$ satisfies the MA condition. Then $\Lambda(A_\Phi(G)^{\ast\ast}) = A_\Phi(G)$ if and only if $\Lambda(A_\Phi(H)^{\ast\ast}) = A_\Phi(H)$ for every open $\sigma-$compact subgroup $H$ of $G.$
    \end{thm}

    \begin{proof}
        Assume that $\Lambda(A_\Phi(G)^{\ast\ast}) = A_\Phi(G).$ Let $H$ be an open $\sigma-$compact subgroup of $G$ and $\Gamma \in \Lambda(A_\Phi(H)^{\ast\ast}).$ Then, by Lemma \ref{ihmh}, we have, $$i_H^{\ast\ast}(\Gamma) \in \Lambda(A_\Phi(G)^{\ast\ast}).$$ So, $i_H^{\ast\ast}(\Gamma) \in A_\Phi(G).$ Now, by elementary calculations, it is easy to verify that the function $i_H^{\ast\ast}(\Gamma)$ vanishes outside $H$ and thus, by the usual identification of $A_\Phi(H)$ in $A_\Phi(G),$ it follows that $\Gamma = i_H^{\ast\ast}(\Gamma)$ belongs to $A_\Phi(H).$ Hence, $\Lambda(A_\Phi(H)^{\ast\ast}) = A_\Phi(H).$

        Conversely, assume that $\Lambda(A_\Phi(H)^{\ast\ast}) = A_\Phi(H)$ for any open $\sigma-$compact subgroup $H$ of $G$ and let $\Gamma \in \Lambda(A_\Phi(G)^{\ast\ast}).$ Then, by Lemma \ref{ihmh}, we have,
        $$m_H^{\ast\ast}(\Gamma) \in \Lambda(A_\Phi(H)^{\ast\ast}).$$ 
        So, $m_H^{\ast\ast}(\Gamma) \in A_\Phi(H).$ Now, by Theorem \ref{restriction}, the restriction of $\Gamma$ (say $\Gamma_r$) to $UCB_\Psi(\widehat{G})$ is in $A_\Phi(G),$ i.e., $\Gamma_r = J(u_\Gamma)$ for some $u_\Gamma \in A_\Phi(G),$ where $J$ is the canonical map from $A_\Phi(G)$ to $A_\Phi(G)^{\ast\ast}.$ Since $G$ is amenable, by \cite[Theorem 3.1]{RLSK2}, the algebra $A_\Phi(G)$ has a bounded approximate identity, say $\{e_\alpha\}$ and further, by Theorem \ref{rightidentity}, the algebra $A_\Phi(G)^{\ast\ast}$ has a right identity, say $\Gamma_e.$ Moreover, it is clear from the proof of \cite[Proposition III.28.7]{BD} that $\Gamma_e$ is the $w^\ast-$ cluster point of $\{J(e_\alpha)\}.$ Thus, for any $T \in PM_\Psi(G),$ since the function $e_\alpha \cdot T$ is in $UCB_\Psi(\widehat{G}),$ we have,
        $$\langle \Gamma, e_\alpha \cdot T \rangle = \langle \Gamma_r, e_\alpha \cdot T \rangle = \langle J(u_\Gamma), e_\alpha \cdot T \rangle = \langle T,u_\Gamma \, e_\alpha \rangle \to \langle T,u_\Gamma \rangle = \langle J(u_\Gamma),T \rangle.$$
        
        Furthermore, since $\{J(e_\alpha)\}$ converges in $w^\ast-$ topology to $\Gamma_e$ and $\Gamma \in \Lambda(A_\Phi(G)^{\ast\ast}),$ it follows that,
        $$\langle \Gamma, e_\alpha \cdot T \rangle = \langle \Gamma, J(e_\alpha) \odot T \rangle = \langle \Gamma \, \square \, J(e_\alpha), T \rangle \to \langle \Gamma \, \square \, \Gamma_e, T \rangle = \langle \Gamma,T \rangle.$$
        This implies that $\Gamma = J(u_\Gamma),$ i.e., $\Gamma$ is in $A_\Phi(G).$ Hence, $\Lambda(A_\Phi(G)^{\ast\ast}) = A_\Phi(G).$
    \end{proof}

    In the subsequent results, we study the topological center of the Banach algebra $UCB_\Psi(\widehat{G})^\ast.$ From \cite[Corollary 3.4]{AD25}, we recall that the space $W_\Phi(G)$ is continuously embedded in $UCB_\Psi(\widehat{G})^\ast,$ under the assumption that the group $G$ is amenable. Further, it is clear from the identification in \cite[Theorem 3.2 and Corollary 3.4]{AD25} that for any $w \in W_\Phi(G),$ the element $\nu(w)$ in $UCB_\Psi(\widehat{G})^\ast$ is given by
    $$\langle \nu(w),u\cdot T\rangle = \langle T,w \, u\rangle,$$
    where $u \in A_\Phi(G)$ and $T \in PM_\Psi(G).$

    We now present the following important lemma.

    \begin{lem}\label{contained}
        For any amenable group $G,$ the space $W_\Phi(G)$ is contained in the topological center of $UCB_\Psi(\widehat{G})^\ast,$ i.e.,
        $$W_\Phi(G) \subseteq \Lambda(UCB_\Psi(\widehat{G})^\ast).$$
    \end{lem}

    \begin{proof}
        To prove this lemma, it is enough to show that the space $W_\Phi(G)$ is contained in the center of $UCB_\Psi(\widehat{G})^\ast$ as in \cite[Proposition 4.5]{LL}. Now, let $w \in W_\Phi(G)$ and $\Gamma \in UCB_\Psi(\widehat{G})^\ast.$ Then, for any $u \in A_\Phi(G)$ and $T \in PM_\Psi(G),$ we have,
        \begin{align*}
            \langle \nu(w) \, \square \, \Gamma,u \cdot T \rangle &= \langle \nu(w),\Gamma \odot (u \cdot T)\rangle\\ &= \langle \nu(w),u \cdot (\Gamma \odot T)\rangle = \langle \Gamma \odot T, w \, u \rangle = \langle \Gamma, (w \, u) \cdot T \rangle.
        \end{align*}
        On the other hand,
        $$\langle \Gamma \, \square \, \nu(w),u \cdot T \rangle = \langle \Gamma, \nu(w) \odot (u \cdot T)\rangle = \langle \Gamma, w \cdot (u \cdot T) \rangle = \langle \Gamma, (w \, u) \cdot T \rangle.$$
        Thus, 
        $$\nu(w) \, \square \, \Gamma = \Gamma \, \square \, \nu(w),$$
        for all $w \in W_\Phi(G)$ and $\Gamma \in UCB_\Psi(\widehat{G})^\ast.$ Hence, the result follows.
    \end{proof}

    The final theorem of the article deals with the topological center of the algebra $UCB_\Psi(\widehat{G})^\ast$ and is analogous to Theorem \ref{topo}. The proof of the theorem follows verbatim to \cite[Corollary 4.6 (ii)]{miao2009} and hence, we omit it. The main ingredients of the proof are Lemma \ref{contained} and the following two facts:
    \begin{enumerate}
        \item Every closed subgroup of an amenable group is amenable.
        
        \item For an amenable group $G,$ the space $PF_\Psi(G),$ which is the predual of $W_\Phi(G),$ is contained in $UCB_\Psi(\widehat{G}).$ This follows from \cite[Corollary 4.13 (i)]{RLSK3}, since $G$ is amenable.
    \end{enumerate}
    Here is the exact statement of the theorem.
    \begin{thm}\label{lastthm}
        Let $G$ be an amenable group. Then $\Lambda(UCB_\Psi(\widehat{G})^\ast) = W_\Phi(G)$ if and only if $\Lambda(UCB_\Psi(\widehat{H})^\ast) = W_\Phi(H)$ for every open $\sigma-$compact subgroup $H$ of $G.$
    \end{thm}

    
    \section*{Acknowledgement}
    The first author would like to thank the Indian Institute of Technology Delhi for providing the Institute Assistantship.

    \section*{Data Availability} 
    Data sharing does not apply to this article as no datasets were generated or analysed during the current study.

    \section*{Competing Interests}
    The authors declare that they have no competing interests.     
	\bibliographystyle{acm}
	\bibliography{Article5}
	
\end{document}